\documentclass{elsarticle}

\usepackage[margin=1in]{geometry}
\usepackage{lineno,hyperref}
\usepackage{amsmath, amssymb}
\usepackage{float}
\usepackage[]{algorithm}
\usepackage{algpseudocode}
\usepackage{graphicx}
\usepackage{appendix}
\usepackage{stmaryrd}
\usepackage{url}

\journal{arXiv (under journal review)}

\usepackage{natbib}
\biboptions{numbers,sort&compress}

\begin{document}
\bibliographystyle{elsarticle-num}
\begin{frontmatter}

\title{Practical Optimal Control of a Wave-Energy Converter in Regular Wave Environments}

\author[ise]{Mertcan Yetkin\corref{mycorrespondingauthor}}
\ead{mey316@lehigh.edu}
\author[mech]{Sudharsan Kalidoss}
\ead{suk515@lehigh.edu}
\author[ise]{Frank E. Curtis}
\ead{frank.e.curtis@lehigh.edu}
\author[ise]{Lawrence V. Snyder}
\ead{lvs2@lehigh.edu}
\author[mech]{Arindam Banerjee}
\ead{arb612@lehigh.edu}
\address[ise]{Department of Industrial and Systems Engineering, Lehigh University, Bethlehem PA, USA}
\address[mech]{Department of Mechanical Engineering and Mechanics, Lehigh University, Bethlehem PA, USA}


\cortext[mycorrespondingauthor]{Corresponding author}

\begin{abstract}
A generic formulation for the optimal control of a single wave-energy converter (WEC) is proposed.  The formulation involves hard and soft constraints on the motion of the WEC to promote reduced damage and fatigue to the device during operation. Most of the WEC control literature ignores the cost of the control and could therefore result in generating less power than expected, or even negative power. Therefore, to ensure actual power gains in practice, we incorporate a penalty term in the objective function to approximate the cost of applying the control force.  A discretization of the resulting optimal control problem is a quadratic optimization problem that can be solved efficiently using state-of-the-art solvers.  Using hydrodynamic coefficients estimated by simulations made in WEC-Sim, numerical illustrations are provided of the trade-off between careful operation of the device and power generated.  Finally, a demonstration of the real-time use of the approach is provided.

\end{abstract}

\begin{keyword}
Wave energy conversion, optimal control, mathematical optimization, constraints, long-term effectiveness, safe operation
\end{keyword}

\end{frontmatter}


\section{Introduction}
\label{intro}

To satisfy growing electricity demands in a sustainable manner, the world needs to rely more on renewable energy resources.  In this context, harnessing the energy from ocean waves can be considered a valuable addition to current renewable energy technology.  The authors in \cite{pecher2017handbook} report that 20\% of the world's total energy demand can be supplied from ocean waves.  Although the early recognized work in wave energy conversion dates back 40 years \cite{budal1977theory,budar1975resonant}, its popularity has increased in the last 10-15 years \cite{antonio2010wave}.  Wave energy conversion technology mainly depends on two elements: design of the wave energy converters (WECs) and optimal control methods for efficient energy extraction.  This paper focuses on the latter.

Researchers are interested in developing novel methods to increase the efficiency of energy extraction from waves using several approaches, including mathematical optimization and optimal control theory \cite{hals2011comparison, ringwood2014energy}.  The work in this paper falls under this category of research.  There are numerous WEC designs employing different conversion technologies \cite{antonio2010wave, drew2009review}
, and energy extraction techniques are inherently dependent on the design.  For our purposes, we focus our scope on point absorbers under heaving motion, due to their potential performance and general applicability.  Specifically, the two-body point absorber described by Muliawan et al. \cite{muliawan2013analysis} is considered throughout this study.  That said, we speculate that the ideas proposed in this paper can be extended by researchers working on optimal control strategies for different types of devices.

\begin{figure}
\centering
\includegraphics[width=6cm]{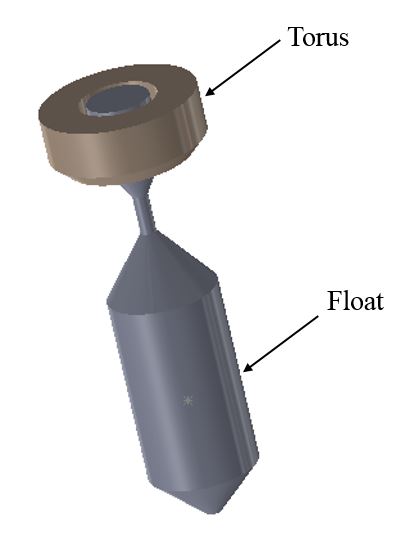}
\caption{Schematic of floating two-body point absorber}
\label{SRPA}
\end{figure}

A schematic of a floating two-body point absorber is shown in Figure \ref{SRPA}.  The float is a slender body running deep in the water while the torus, with a shallower body shape, floats on the surface.  The float and torus are wave-following devices and oscillate at different frequencies with the wave.  The relative motion between the two in the heave (vertical) direction is exploited for energy conversion.

Optimizing energy extraction of WECs can be classified into two essential components, commonly described as frequency-domain control and time-domain control \cite{falnes2002ocean}.  Often, time-domain control approaches are desirable since they present the opportunity to track and adjust the state variables including, but not limited to, position, velocity, and acceleration. In the present work, we propose a generic model that can impose bounds on these state variables, which is possible under a time-domain approach.

Optimal control methodologies have been studied in the literature for optimizing WEC operation. Falnes \cite{falnes2002optimum} provides a general review of optimal control methods in wave energy extraction and provides a thorough historical progression. Hals et al. \cite{hals2011constrained} formulate the optimal control problem  by considering bounds on the position and the control force and further employ a model-predictive control (MPC) strategy. These authors also note that optimal energy absorption may be better in irregular waves compared to regular waves having the same energy-period and wave-power level. Wang et al. \cite{wang2015constrained} formulate a constrained optimal control problem by imposing limitations on the motion and force while having the control force as a spring-damper controller with a time-dependent damper term. Further, Zou et al. \cite{zou2017optimal} consider a heaving WEC and maximize the energy extraction with both unconstrained and constrained optimal control methodologies. For efficient real-time implementation on real sea data, Bacelli et al. \cite{bacelli2011control} consider a formulation with bounds on the displacement and approximate the states and the force by a linear combination of basis functions.

Model predictive control (MPC), which operates on a moving horizon, presents a framework to adjust the problem by leveraging new information and, thus, aids in the implementation of real-time control.  Cretel et al. \cite{cretel2011maximisation} present an MPC-based control method for maximizing the power output, alleviating the drawbacks of current control strategies such as latching control and reactive control. Richter et al. \cite{richter2013nonlinear} argue that a more complex MPC strategy is needed in order to incorporate nonlinear effects such as mooring forces.  To this end, they propose a nonlinear model predictive control (NPMC) approach as an alternative to the current MPC techniques.  Nielsen et al. \cite{nielsen2017optimizing} adopt an MPC framework considering not only the power generation but also the losses in the mechanical-to-electrical conversion of the power. Arguing that the fluctuation in the control strategy causes fatigue in some of the components of WEC, a comparison is made between the power generated and the fatigue.

Throughout these related works on optimal control and/or MPC, little attention has been given to some issues that arise when employing control strategies in practice, such as the cost of applying a control force, inefficiencies in the power conversion system, and long-term damage to the device. In order to make wave energy conversion a viable part of the renewable energy portfolio, more work needs to be done to address such issues. To this end, the primary objectives of this work presented in this paper are the following:
\begin{itemize}
    \item We augment an optimal control formulation for a WEC by including hard and soft constraints on the device.  The inclusion of the constraints of this type can aid in reducing possible damage to the device including short-term stress and long-term fatigue over the components. A more detailed explanation can be found in Section~\ref{safety_constraints}.
    \item We  augment the objective function with a  term designed to account for energy expended in applying the control force.  This cost of the control force plays a crucial role.  Often ignored in other formulations, our experiments demonstrate that, without accounting for this cost, the net energy could be much less than the model predicts---or even negative.
    \item We demonstrate that our formulation can be solved quickly enough to execute in real time in a regular wave with experiments over a rolling horizon.
\end{itemize}

The remainder of the paper is organized as follows. Section \ref{hydro} presents the hydrodynamic model for the two-body point absorber.  Section \ref{formulation} describes an optimal control formulation of the problem and further details. Section \ref{real-time} discusses a methodology for real-time implementation of the control strategy.  Numerical results are presented in Section \ref{numerical_results}.  Concluding remarks are provided in Section \ref{conclusion}. 


\section{Hydrodynamic model of a point absorber}
\label{hydro}


\subsection{Governing equation}
\label{governing_equation}
The governing equation of motion of a single point absorber in the time domain, also known as the Cummings equation \cite{cummins1962impulse}, can be written as
\begin{equation}
\label{governing}
    (M + m)\ddot{z}(t) + \int_0^t K(t- \tau)\dot{z}(\tau)d\tau = F_\text{ex}(t) + F_\text{h}(t) + F_\text{PTO}(t) + F_\text{es}(t),
\end{equation}
where $ \ddot{z}(t)$ and  $\dot{z}(t)$ respectively represent the acceleration $(\text{m}/\text{s}^{2})$ and velocity $(\text{m}/\text{s})$ of the body at time $t$ in the heave mode, $M$ (kg) is the mass of the body, $m$ (kg) is the added mass (representing the mass of fluid displaced by the body during operation), and $K(\cdot)$ $(\text{kg}/\text{s}^{2})$ is the memory function of the radiation force that represents the hydrodynamic damping of the system. The forces considered for the current analysis include the excitation force on the body, $F_\text{ex}$ (N); the hydrostatic force, $F_\text{h}$ (N); the power take-off (PTO) force, $F_\text{PTO}$ (N); and, the end stop force $F_\text{es}$ (N) ensuring some stopping conditions. We assume linear wave theory and ignore viscous drag on the system \cite{kalidoss2019performance}. 

The excitation force can be written as
\begin{equation}
    \label{excitation2}
    F_\text{ex}(t) = \Re \left( \sum_i \sqrt{2 S(f_i) \Delta f}  \bar{F}_\text{ex}(f_i) e^{i \left( 2 \pi f_i t + \phi_i \right)} \right)
\end{equation}
where $f_i$ ($1/\text{s}$) are wave frequencies, $\phi_i$ (rad) are random phases, $\bar{F}_\text{ex}$ ($\text{N}/\text{m}$) represent complex vectors of wave excitation force per meter of wave amplitude in the frequency domain, $\Delta f$ ($1/\text{s}$) is an appropriate frequency step, and $S(f)$ ($\text{m}^{2}/\text{Hz}$) represents the spectrum of the incident wave field \cite{babarit2012numerical}.  For single regular wave conditions investigated in this paper, the excitation force $F_\text{ex}$ is linear in terms of $w(t)$.

The hydrostatic force can be written as
    \begin{equation}
    \label{hydrostatic}
        F_\text{h}(t) = \rho g S_w z(t),
    \end{equation}
    where $\rho$ ($\text{kg}/\text{m}^3$) is the density of the water, $g$ ($\text{m}/\text{s}^2$) is the acceleration due to gravity and $S_w$ ($\text{m}^2$) is the cross section area of the WEC at the water plane. Since the parameters in equation \eqref{hydrostatic} are constant, it can be rewritten as
    \begin{equation}
        \label{hydrostatic_mod}
        F_h(t) = K_h z(t).
    \end{equation}

Due to the mechanical complexity of the PTO system, it is common to limit operating stroke length of the WEC by introducing physical bounds, known as end stops \cite{hong2016impact}. Numerically, the end stops are modeled as a restoring force with large restoring coefficient \cite{muliawan2013analysis}. The end stop force activates only when the relative motion between the WEC bodies exceeds the limits; for the current simulations, a limit of $-3$ m to $3$ m is used. The end stop force acting on the WEC is given as
    \begin{equation}
        \label{end stop}
        F_{es} = -K_{es}\left( (z(t)+z_{es})u(-z_{es}-z(t)) + (z(t)-z_{es})u(z(t)-z_{es}) \right),
    \end{equation}
    where $K_{es}$ ($\text{N}/\text{m}$) is the restoring coefficient, $z_{es}$ (m) is the stroke limit and $u$ is the step function. In most operating conditions, the stroke length is within the limit.
    
    In the present analysis, the point absorber is attached to the seabed with a long mooring line so that the mooring force acting on the system is negligible. The power take-off (PTO) force is ideal in the sense that it can exert any desired force within some limits.


\subsection{Hydrodynamic model used in the simulation}
 We estimate the system coefficients for the control system from a numerical simulation of the point absorber. In our numerical simulation, the hydrodynamic coefficients, such as added mass, hydrodynamic damping coefficient and excitation forces, are calculated in the frequency domain using WAMIT \cite{lee2013wamit}, which is a potential flow solver based on a boundary element method (BEM). The linear wave theory assumption is adopted in the BEM solver, i.e., the wave amplitude is assumed to be very small compared with the wavelength. Therefore, the body-boundary conditions and free surface are linearized, and the gradient of velocity potential, $\Phi$, satisfies the Laplace equation ($\nabla^2\Phi = 0$) in the fluid domain. The added-mass and damping coefficients calculated using velocity potential are given by 
\begin{equation}
    \label{potential_eqn}
    m_{ij} -\frac{i}{\omega} \zeta_{ij} = \rho\iint_{S_b} n_i\phi_jdS,
\end{equation}
where $m_{ij}$ is the added-mass coefficient and $\zeta_{ij}$ is the hydrodynamic damping coefficient.

The excitation force ($X_i$) on the WEC due to the incident wave potential ($\phi_0$) is given by 
\begin{equation}
    \label{ex_force}
    X_i = -i\omega\rho\iint_{S_b}\left( n_i\phi_0 - \phi_i\frac{\partial\phi_0}{\partial n}\right)dS.
\end{equation}

The multi-body dynamics of the WEC is modeled using WEC-Sim \cite{ruehl2014preliminary}, which is a MATLAB-based, open-source code. In WEC-Sim, the complex interaction between the incident waves, device motion, and power take-off mechanism is modeled in the time domain. The power performance and device motion of the WEC are simulated using the radiation and diffraction method. In the present analysis, the hydrodynamic bodies (float and torus) of the WEC are surface meshed and constrained to operate in one degree of freedom (heave, in the present case). With the hydrodynamic coefficients calculated from WAMIT, WEC-Sim calculates the time-domain motion and power performance of the WEC.
In Section \ref{estimation}, we describe how we estimated the  system coefficients, including the hydrodynamic system, control, and wave coefficients.

\section{Mathematical formulation of the optimal control problem}
\label{formulation}
In this section, we develop an optimal control formulation.  We first simplify the governing equation using state variables, and then introduce a basic formulation with additional components designed to promote desirable solution properties.
Our complete formulation is given in Section \ref{sec.complete_formulation}. 

By using the simplifications introduced in Section \ref{governing_equation}, we have the refined version of the equation \eqref{governing} of the WEC with state variables as
\begin{equation}
\label{governing_linear}
    (M + m)\ddot{z}(t) + K \dot{z}(t) = K_\text{ex}w(t) + K_\text{h} z(t) + F_\text{PTO},
\end{equation}
where $K = K(\cdot)$ and $K_{ex} w(t) = F_{ex}(t)$ for some $K_{ex}$ with only a regular wave as input.  We shall simplify further by approximating the acceleration term involving $\ddot{z}(t)$.  Not that, for a sinusoidal wave, say defined by $\kappa(t) = \sin(at)$, the position and acceleration have perfect negative correlation due to the fact that $\ddot\kappa(t) = -a^2\sin(at)$.  
Thus, we can further simplify \eqref{governing_linear} to obtain
\begin{equation}
\label{governing_linear_without_acc}
    K \dot{z}(t) +(- K_\text{h} + \xi)z(t) = K_\text{ex}w(t) + F_\text{PTO},
\end{equation}
where $\xi z(t)$ approximates the acceleration term from \eqref{governing_linear}.

In the remainder of this section, we develop our proposed optimal control formulation.  We begin with a basic formulation, then add components designed to promote desirable solution properties.  Our complete formulation is given in Section \ref{sec.complete_formulation}.

\subsection{Basic optimal control problem}
\label{basic_formulation}
The objective of WEC control is the maximization of the total energy extracted throughout a given horizon $T$, which can be simply described as
\begin{equation}
    \label{continuous_max_objective}
    \max_{u,\dot{z},z} -\int_0^T u(t)\dot{z}(t)dt,
\end{equation}
where $u(t)$ is the control force at time $t$.
In the context of the governing equation described in the previous section, $u(t)$ represents the PTO force at time $t$. Based on  \eqref{governing_linear_without_acc}, the control force is considered to have a linear relationship with the states of the system.

Applying a forward Euler discretization to the system dynamics in \eqref{governing_linear_without_acc}, we have the following discrete-time optimal control problem:
\begin{equation}
     \label{basic_obj} 
     \max_{u,\dot{z},z} -\sum_{k=1}^{N}{\frac {[u(t_{k-1})\dot{z}(t_{k-1})+u(t_{k})\dot{z}(t_{k})]}{2}}\Delta t_{k}
 \end{equation}
 subject to:
\begin{align}
&\text{a given initial point $(\dot{z}(0),z(0))$}; \label{initial} \\
&\left[ \begin{array}{c}
    \label{discrete_system} \dot{z}(t_{k+1})       \\
        z(t_{k+1})  
    \end{array} \right ] 
    = A \left [ \begin{array}{c}
    \dot{z}(t_k)       \\
        z(t_k)  
    \end{array} \right ] + b u(t_k) + c w(t_k)  &\text{for } k \in \{0,1,\dots,N\}; \\
 \label{hard_bound} &-\gamma \leq u(t_k) \leq \gamma &\text{for } k \in \{0,1,\dots,N\},
\end{align}
where $\Delta t_k$ is the discretization interval for the $k$th time period, $k \in \{0,1,\dots,N\}$. In the discretized system dynamics represented by \eqref{discrete_system}, $A \in \mathbb{R}^{2 \times 2}$ is the coefficient matrix for the discretized system involving the velocity and the position of the device, $b \in \mathbb{R}^2$ is the coefficient vector for the control force $u(t)$, and $c \in \mathbb{R}^2$ is the coefficient vector for the wave position $w(t)$. Equation \eqref{initial} enforces the initial condition of the system given an initial point, whereas \eqref{hard_bound} limits the control force applied into the system by $\gamma$, denoting the maximum physical limit on the force, in Newtons. We refer to this constraint as a ``hard bound.''

Observe that $\Delta t_k$ can be different for different $k$, based on the specific discretization method used. In the formulation above, the objective in  \eqref{basic_obj} is approximated via the trapezoidal method. If an equidistant discretization method is considered, the objective in \eqref{basic_obj} would be given by
\begin{equation}
    \label{simplified_basic_obj}
    \max_{u,\dot{z},z} -\sum_{k=0}^{N} u(t_{k})\dot{z}(t_{k}),
\end{equation}
which is common in the literature; see, e.g., \cite{hals2011constrained,cretel2011maximisation}.

\subsection{Constraints to promote reduced damage to the device}
\label{safety_constraints}

The optimal solution to the formulation above may involve control forces that are undesirable from the perspective of the safety and longevity of the device. 
One concern is that the force might alternate quickly between a large positive quantity and a large negative quantity, creating a hammering effect that can cause fatigue over the long term. Another concern is that the force might stay at a large positive or a large negative quantity for a long time and apply constant pressure to the mechanical components, especially to the joint parts of the device. This, too, may result in undue wear-and-tear on the device. 

A solution to these potential causes of damage is to impose a ``soft bound'' on the control force, wherein one penalizes the objective function  for any control value that exceeds a certain fraction of the hard bound constraint \eqref{hard_bound}.  In particular, we incorporate into the objective the term
\begin{equation}
\label{soft_bound}
    - \sum_{k=1}^N \frac{\rho \left (|u(t_{k-1})| - \gamma \eta \right )^+ + \rho \left (|u(t_{k})| - \gamma \eta \right )^+ }{2} \Delta t_k,
\end{equation}
where $\rho \geq 0$ is the per unit penalty for exceeding the soft bound, $(\cdot)^+ = \max\{0,\cdot\}$, and $\eta \in (0,1]$ is the fraction of the hard bound at which to impose the soft bound.

Finally, it is very common to have a physical limit on the displacement of the torus relative to the float. As is done in several works \cite{bacelli2011control, hals2011constrained,richter2013nonlinear,wang2015constrained}, we consider a bound on the maximum displacement of the device (which is in fact the relative displacement between the torus and the float), given as
\begin{equation}
\label{displacement_bound}
    - \delta \leq z(t_k) \leq \delta \qquad \text{for } k \in \{0,1,\dots,N+1\}, 
\end{equation}
where $\delta$ is the constant of displacement in meters. Another interpretation of \eqref{displacement_bound} is that it prevents the torus from going underwater completely, or rising above the sea level and slamming onto the sea surface.

\subsection{Cost of applying the control force}
\label{cost_of_control}

Many WEC control models in the literature do not account for the cost or energy required to apply the control. In other words, they maximize the energy captured by the WEC using formulations similar to \eqref{basic_obj}--\eqref{hard_bound}, but they ignore the energy that must be expended to exert this control force. As a result, the net energy captured by the device may be much less than predicted by the model---it may even be negative. 

\begin{figure}
\centering
\includegraphics[width=10cm]{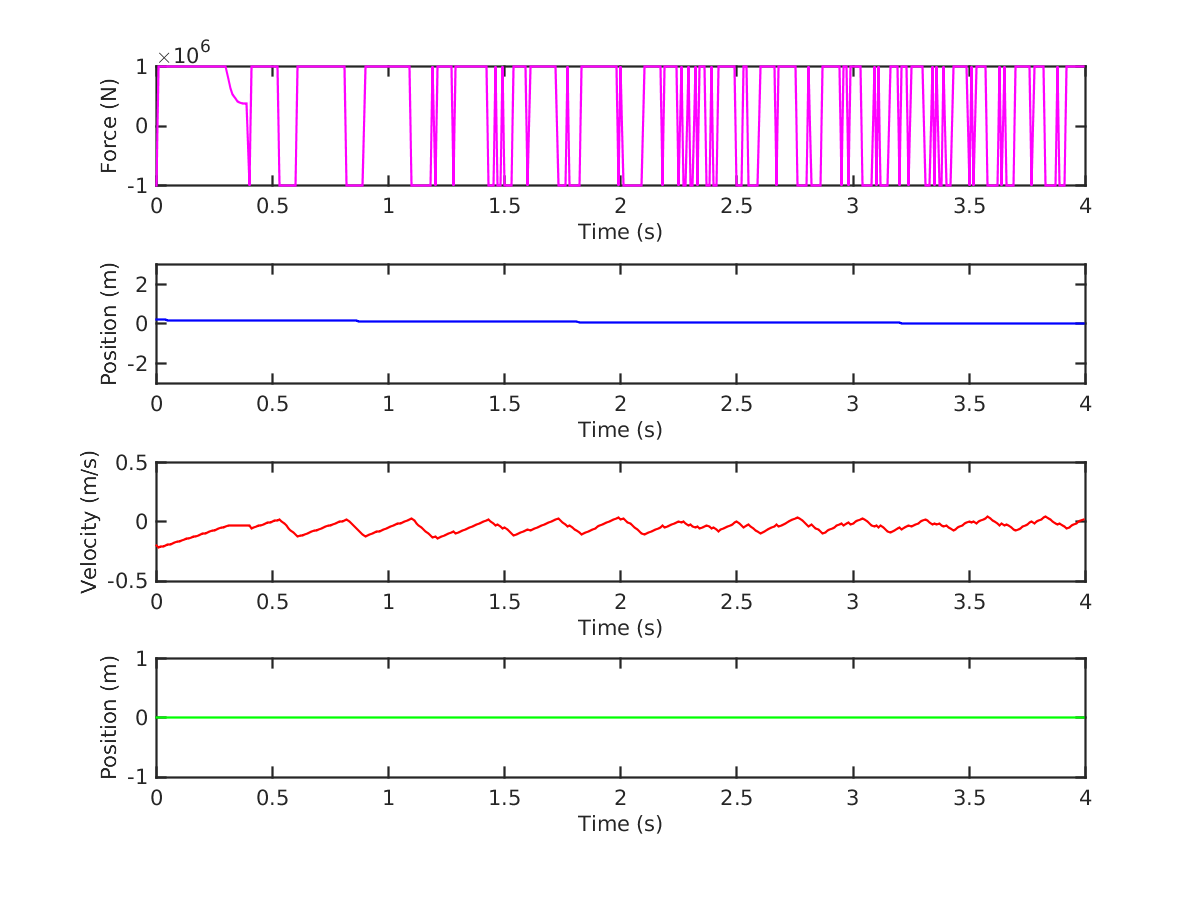}
\caption{Optimal solution of an instance without incorporating the cost of control force. Subplots represent trajectories of control ($u$), position ($z$), velocity ($\dot{z}$) and wave elevation ($w$), respectively.}
\label{zerowave}
\end{figure}

As a simple example, suppose that there is no incident wave, i.e., $F_\text{ex}(t)\equiv 0$. For certain values of the remaining parameters, and with initial conditions, the optimal solution to \eqref{basic_obj}--\eqref{hard_bound} (plotted in Figure~\ref{zerowave}) has non-zero control force and a strictly positive optimal objective function value $\approx 10^5$ (i.e., a strictly positive energy extraction). This suggests that we can extract energy even if there is no wave. Of course, this is erroneous, and the discrepancy comes from the model's failure to account for the energy expended to exert the control force. In this simple example, the motion of the WEC comes solely from the control force (since there is no incident wave), and the energy required to exert that force must exceed the energy extracted from the motion. Hence, this solution involves a net loss of energy, even though the model suggests a gain. Any optimization model that does not incorporate the cost of control will suffer from this same sort of error.

Although the ``cost'' of the control force has not been included explicitly in WEC control models, one could argue that it is accounted for implicitly in some control strategies. For example, relatively simpler control methods such as latching control \cite{babarit2006optimal} or phase control \cite{antonio2008phase} restrict the control policy to have a certain form, rather than optimizing over all possible controllers; since these control methods are common, one might infer that they do not incur excessive control costs. However, we argue that these control methods are too simplistic, since they do not promote device safety and offer less flexibility to be extended to irregular wave environments. Therefore, advanced control models such as the one proposed here are more promising. For such models, including the cost of control as an explicit term in the model is essential.


Of course, it may be difficult to derive an explicit expression for the control cost. Such expressions may be highly complex, and may depend on the specific power conversion system. Therefore, as a proxy, we assume that the instantaneous cost of applying a force $u$ has the functional form $\mathcal{O}(u^2)$, i.e., that it equals $\lambda_1u^2$ for a constant $\lambda_1$. We chose a form by calibrating various nonlinear models using data from \cite{klamo2007effects}, which considers magnetic damping as a controller for a specific system. More details on our derivation of this functional form are given in Appendix \ref{derivation}. Of course, for systems that use different PTO systems and control strategies than those in \cite{klamo2007effects}, the control cost may have a form that is different from $\mathcal{O}(u^2)$. For such systems, one can use the approach we outline in  Appendix \ref{derivation}, or another approach, to estimate the functional form, and then insert that form into our optimization model. That is, our framework is flexible enough to handle any form of the instantaneous control cost, and our contribution lies in the modeling and solution approach, rather than in the specific form used for this cost.


\subsection{Ensuring smoothness of the control force}
\label{smoothing-control}
The optimal control strategy resulting from the basic formulation in Section~\ref{basic_formulation} is likely to be non-smooth since any arbitrary control policy satisfying the bounds is allowed. For example, in bang-bang control policies, the control is always at its bounds.  Some works in the literature, e.g.,  \cite{li2012wave}, try to exploit this type of structure and develop  efficient procedures for finding the optimal solution within the family of bang-bang policies.  However, control policies, like bang-bang, that have sharp jumps in the control might not be desirable for two main reasons: 
\begin{enumerate}
    \item It may not be possible to apply a non-smooth control policy, since most  power take-off systems are best suited for smooth control forces.
    \item Non-smooth control would create additional wear and tear on the components of the WEC.
\end{enumerate}
Thus, we would like the control sequence to be smooth in our formulation. One straightforward way to encourage smoothness is to add the term
\begin{equation*}
   - \lambda_2 (\Delta u(t_k))^2
\end{equation*}
to the objective function, where $\lambda_2$ is a parameter to adjust the importance of smoothness. By approximating $\Delta u(t_k)$ from the differences in $u(t_k)$ and applying the trapezoidal rule, one obtains the term
\begin{equation}
    \label{smooth_control_obj}
    - \frac{\lambda_2 \left( (u(t_{k}) - u(t_{k-1}))^2 + (u(t_{k-1}) - u(t_{k-2}))^2 \right)}{2 \Delta t_k }.
\end{equation}
From another perspective, this term can also be considered as the cost of changing the control force, with a cost coefficient of $\lambda_2$.

\subsection{Complete formulation of the problem}\label{sec.complete_formulation}
Combining all the elements discussed in Sections \ref{basic_formulation}--\ref{smoothing-control}, we have the complete formulation of the discrete-time optimal control problem:
\begin{equation}
     \label{complete_obj} 
     \begin{split}
         \max_{u,\dot{z},z,\alpha}  &\sum_{k=2}^{N} {\frac {-[u(t_{k-1})\dot{z}(t_{k-1})+u(t_{k})\dot{z}(t_{k})] - \rho(\alpha_{k-1} + \alpha_k) }{2}}\Delta t_{k} \\ &-  \frac{\lambda_1(u(t_{k-1})^2 + u(t_{k})^2)}{2}\Delta t_k \\
         & - \frac{\lambda_2 \left( (u(t_{k}) - u(t_{k-1}))^2 + (u(t_{k-1}) - u(t_{k-2}))^2 \right)}{2 \Delta t_k }
     \end{split}
 \end{equation}
 subject to
\begin{align}
&\text{constraints \eqref{initial}, \eqref{discrete_system}, \eqref{hard_bound}, \eqref{displacement_bound}} \nonumber \\
    \label{penalty} &-\alpha_{k} - \gamma \eta \leq u(t_k) \leq \alpha_{k} + \gamma \eta  &\text{for } k \in \{0,1,\dots,N\}\\
    \label{nonneg} & \alpha_k \geq 0 &\text{for } k \in \{0,1,\dots,N\}.
\end{align}
The objective function \eqref{complete_obj} maximizes the total energy extraction minus the cost of applying the control force and penalties for violations of the soft bounds on the control limits.  Constraints \eqref{penalty} and \eqref{nonneg} are used in conjunction with the second term in the objective as an equivalent linear representation of the soft bound condition \eqref{soft_bound}.  Specifically, the auxiliary variables $\alpha_k$ describe the amount by which the force exceeds the soft bound, which are penalized in the objective with a constant per unit cost $\rho$.

We make the following remarks about this complete model.
\begin{itemize}
    \item It is assumed that the unit cost violation term $\rho$ is a constant. In order to achieve different bound structures, one can use a function $\rho(\cdot)$ instead of specifying only one cost coefficient $\rho$, where the function depends on the amount of violation. 
    \item Power smoothness of the operation can be handled internally in the mechanical/electrical power take-off component \cite{Torres2015power, Hansen_2013}.
    \item The same mechanical/electrical component is used for both taking power from the system and feeding power into the system to apply the control.
    \item Since the parameter $\lambda_1$ gives the cost of applying the control force, it also captures the conversion inefficiency of the system. For instance, if we assume a conversion efficiency ratio $0<\psi<1$, then we have a system that can only capture a proportion $\psi$ of the total energy injected into the system by the controller, and we can extract only a proportion $\psi$ of the total energy harvested. In this case, the inefficiency in the conversion system can be incorporated into $\lambda_1$ by properly adjusting this coefficient.
\end{itemize}

\subsection{Estimation of the coefficients}
\label{estimation}

We simulate the wave--structure interaction of the WEC using WEC-Sim. The time-domain simulation of a wave--structure interaction is generated by the radiation and diffraction method. In this method, the hydrodynamic coefficients are calculated with a frequency-domain boundary element method and the system dynamics is solved in the time domain. We then estimate the coefficients $A$, $b$, and $c$ in the linear equation \eqref{discrete_system} by simulating the device under different wave conditions and control policies. Specifically, we use the following three-step procedure. The estimations in each step are done using least-squares regression.
\begin{enumerate}
    \item First, we consider a free-decay experiment in which the device is given an initial condition with no wave input. From this experiment, we estimate the hydrodynamic coefficient matrix $A$.
    \item Second, we consider a free-floating experiment in which the device is floating with some wave input.  Using the previously estimated coefficient matrix~$A$, we estimate the wave coefficient vector~$c$ in this step.
    \item Finally, we conduct an experiment with the same input wave and some control policy. Using the previously estimated coefficients $A$ and $c$, we estimate the control coefficient vector $b$. A caveat here is that we cannot simply use the control policy used in the simulation to estimate the coefficients, since the system representation in \eqref{discrete_system} is only an approximation.
    That said, we know that in regular waves, the control policy has a sinusoidal form with the same period as the wave input.  Thus, to estimate $b$, we consider sinusoidal control policies with a fixed amplitude (as the force limit), with the same period as the wave input, and with all possible shifts within the discretization level; we then choose the coefficients that resulted in the least error in the regression.
\end{enumerate}
Note that all the coefficients are estimated independently to ensure that they accurately represent the system. An advantage of this procedure is that the hydrodynamic coefficient matrix $A$ is estimated only once since the device remains the same.  However, if a different wave input is considered, one needs to start over from step 2. If the new wave input has the same period, it is straightforward to scale the wave coefficients $c$ accordingly. However, it is not so clear how to adjust the coefficients when a wave input with a different period is considered. This problem may be a possible future direction.

Typical wave conditions for the east coast of the United States have an average wave period of 4--6 seconds and significant wave height around 6 meters \cite{noahdata}. We use these numbers to guide our choice of wave inputs in our estimation, which are summarized in Table \ref{wave-inputs}.  Note that these conditions are off-resonance for the device $(\text{wave period} \neq 10s)$. 

\begin{table}
\centering
\begin{tabular}{ccc} \hline
       & \textbf{Significant Wave Height (m)} & \textbf{Period (s)} \\ \hline
Wave 1 & 6                                    & 4                   \\
Wave 2 & 6                                    & 5                  \\
Wave 3 & 6                                    & 6   \\\hline              
\end{tabular}
\caption{Wave inputs considered in the estimation}
\label{wave-inputs}
\end{table}

\subsection{Coefficients for the system}
\label{coefficients_system}

The estimated coefficients for the system in which the input wave has 6-meter significant wave height and 4-second wave period are:
\begin{equation}
    \label{coefficient_values_h_6_t_4}
    A = \left [ \begin{array}{cc}
        0.9939 & -0.0378 \\
        0.00997 & 0.9998
    \end{array}\right], \quad b = \left [ \begin{array}{c}
         0.0123 \times 10^{-6}  \\
         6.1785 \times 10^{-11}
    \end{array} \right], \quad c = \left [\begin{array}{c}
         0.0045 \\
         2.2480 \times 10^{-5} 
    \end{array} \right].
\end{equation}

\begin{figure}
\centering
\includegraphics[width=10cm]{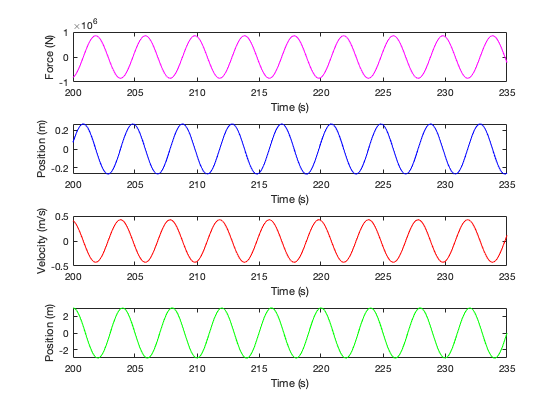}
\caption{Reconstruction of the system trajectory for verification, where the input wave has 6-meter significant wave height and 4-second period. Subplots represent control, position, velocity and wave elevation, respectively.}
\label{iterated-trajectory}
\end{figure}

\begin{figure}
\centering
\includegraphics[width=10cm]{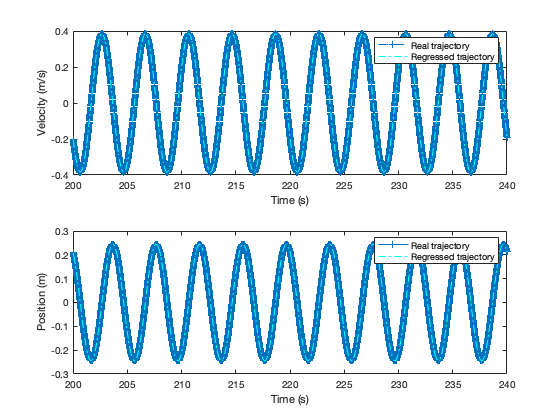}
\caption{Simulation data and regressed trajectories for velocity (first subplot) and position (second subplot) where a control strategy is present and the input wave has 6-meter significant wave height and 4-second period.}
\label{regression-control}
\end{figure}

In order to visually inspect the dynamics of the hydrodynamic system using the coefficients from \eqref{coefficient_values_h_6_t_4}, we rolled forward the equations in \eqref{discrete_system}. Figure \ref{iterated-trajectory} plots the resulting trajectories using using a sinusoidal control policy that respects the bounds on the control variables.  (The control policy used has a similar form to the policy used in the third phase of the estimation). From the figure, it is clear that the system appears to be stable, with trajectories that conform to expectations.  Moreover, by comparing the real and computed trajectories in Figure \ref{regression-control}, it is clear that the estimated coefficients yield dynamics that accurately represent the system.

The coefficients for 6-meter significant wave height and 5-second wave period are:
\begin{equation}
    \label{coefficient_values_h_6_t_5}
    A = \left [ \begin{array}{cc}
        0.9939 & -0.0378 \\
        0.00997 & 0.9998
    \end{array}\right], \quad b = \left [ \begin{array}{c}
         0.0325 \times 10^{-6}  \\
         1.6256 \times 10^{-10}
    \end{array} \right], \quad c = \left [\begin{array}{c}
         0.0142 \\
         7.0887 \times 10^{-5} 
    \end{array} \right].
\end{equation}
Finally, the coefficients for 6-meter, 6-second waves are:
\begin{equation}
    \label{coefficient_values_h_6_t_6}
    A = \left [ \begin{array}{cc}
        0.9939 & -0.0378 \\
        0.00997 & 0.9998
    \end{array}\right], \quad b = \left [ \begin{array}{c}
         0.0429 \times 10^{-6}  \\
         2.1485 \times 10^{-10}
    \end{array} \right], \quad c = \left [\begin{array}{c}
         0.0204 \\
         1.0219 \times 10^{-4} 
    \end{array} \right].
\end{equation}

\section{Receding horizon approach for real-time control}
\label{real-time}

In this section, we discuss our receding-horizon approach for WEC control. The idea is that, rather than solving the optimization problem over a long time horizon, we instead solve reasonably small problems over short time horizons, during which our estimates of the wave inputs and the state of the dynamical system are likely to be accurate. Then, we implement some portion of the resulting control strategy, then make a new forecast of the wave over the next time horizon, and re-solve the model. 

In particular, suppose we wish to solve problems with a horizon length of $T$ time periods, discretized into time steps of length $\Delta t$. Every $T^c$ time periods (called the ``controller update horizon''), we re-solve the problem. For example, if $T=100$ seconds, $\Delta t = 0.1$ second, and $T^c = 5$ seconds, then every 5 seconds, we solve a problem consisting of $T / \Delta t=1000$ discretization points. We implement the control values for the first 5 seconds of that horizon, then update the wave forecast and re-solve. 

The choice of the optimization horizon $T$ is crucial in practice, since one wants it large enough so that the resulting strategy accounts for a longer period so it is not myopic and some part of it can be implemented before the next problem is set up and solved. On the other hand, one wants $T$ to be small enough that the current estimates of the system dynamics and wave are accurate. Finally, we must set the discretization time step $\Delta t$ and the controller update horizon $T^c$ large enough so that the horizon-$T$ problem can be solved within $T^c$ time units. We refer to the model's ability to solve the horizon-$T$ problem within the available time, i.e., within $T^c$, as {\em real-time control.}

\begin{algorithm}
\caption{Receding horizon control strategy}
\label{receding_horizon_alg}
\begin{algorithmic}[1]
\State \textbf{Input:} State at initial time $T_0$, optimization time horizon~$T$, discretization time step~$\Delta t$, and controller update horizon~$T^c$
\State \textbf{Output:} Sequence of control values
\For{$j=0,1,2,\dots$}
\State Set the current time as $T_j$
\State \label{step:predict} Obtain a wave prediction for the next $T$ seconds
\State Solve the constrained optimal control problem over the optimization time horizon using the initial state at $T_j$, the discretization time step $\Delta t$, and the wave prediction over $[T_j,T_j + T]$
\State Implement the control values for the next $T^c$ seconds and collect the state at time $T_j + T^c$
\State Update $T_{j+1} \leftarrow T_j + T^c$
\EndFor
\end{algorithmic}
\end{algorithm}

Algorithm~\ref{receding_horizon_alg} presents our strategy for real-time control.  To handle the general situation of irregular waves, the algorithm includes a step to make short-term wave prediction \cite{fusco2010short}.
For simplicity, however, in the remainder of our study, we consider regular waves.  This obviates the need to perform step \ref{step:predict} and to modify the dynamical system coefficients with each time horizon.

\section{Numerical results}
\label{numerical_results}

In this section, we provide the results of numerical experiments for solving the optimization problem given in Section~\ref{sec.complete_formulation}. Firstly, in Section~\ref{choosing-lambdas}, we discuss a strategy for choosing the cost coefficients $\lambda_1$ and $\lambda_2$, providing numerical demonstrations of experiments used when choosing these coefficients. We also demonstrate the sensitivity of the objective value with respect to $\lambda_1$ and $\lambda_2$ to further emphasize the importance of these cost coefficients. Then, in Section~\ref{choosing-safety-parameters}, we experiment over different values of the parameters $\eta$ and $\rho$ and demonstrate their effect on the difficulty of the problem and the objective value. Lastly, in Section~\ref{receding-horizon-experiments}, we employ the receding horizon control scheme and demonstrate the capabilities of our model in a more realistic setting.

The mathematical models are formulated in the AMPL modelling language \cite{fourer1993ampl} and solved via IPOPT \cite{wachter2006implementation}, which employs a primal--dual interior point method. Furthermore, all the experiments in this section are done considering a wave input with 6-meter significant wave height and 4-second wave period. We also conducted our experiments on instances with other wave inputs, but we did not include the results of those experiments since the main conclusions are similar. 

\subsection{Choosing the cost coefficients $\lambda_1$ and $\lambda_2$}
\label{choosing-lambdas}
The values of the objective function coefficients $\lambda_1$ and $\lambda_2$, which are intended to capture the cost of applying and changing the control force, are difficult to quantify precisely, but it is possible to identify a range of plausible values for them. Note that if we set $\lambda_1$ or $\lambda_2$ to too small a value, then we are saying that the control force can be applied nearly for free; this will result in the control force altering the device motion dramatically, not following the regular period of the waves, as described in Section \ref{cost_of_control}. For sufficiently large values of the coefficients, the control force is so expensive that we will not apply it at all, and the device will simply follow the  wave. With a regular wave input, the velocity of the device should have the same period as the wave; therefore, our suggestion is to choose the smallest values of $\lambda_1$ and $\lambda_2$ that result in the period of the velocity being sufficiently close to that of the wave. 

\begin{figure}
\centering
\includegraphics[width=10cm]{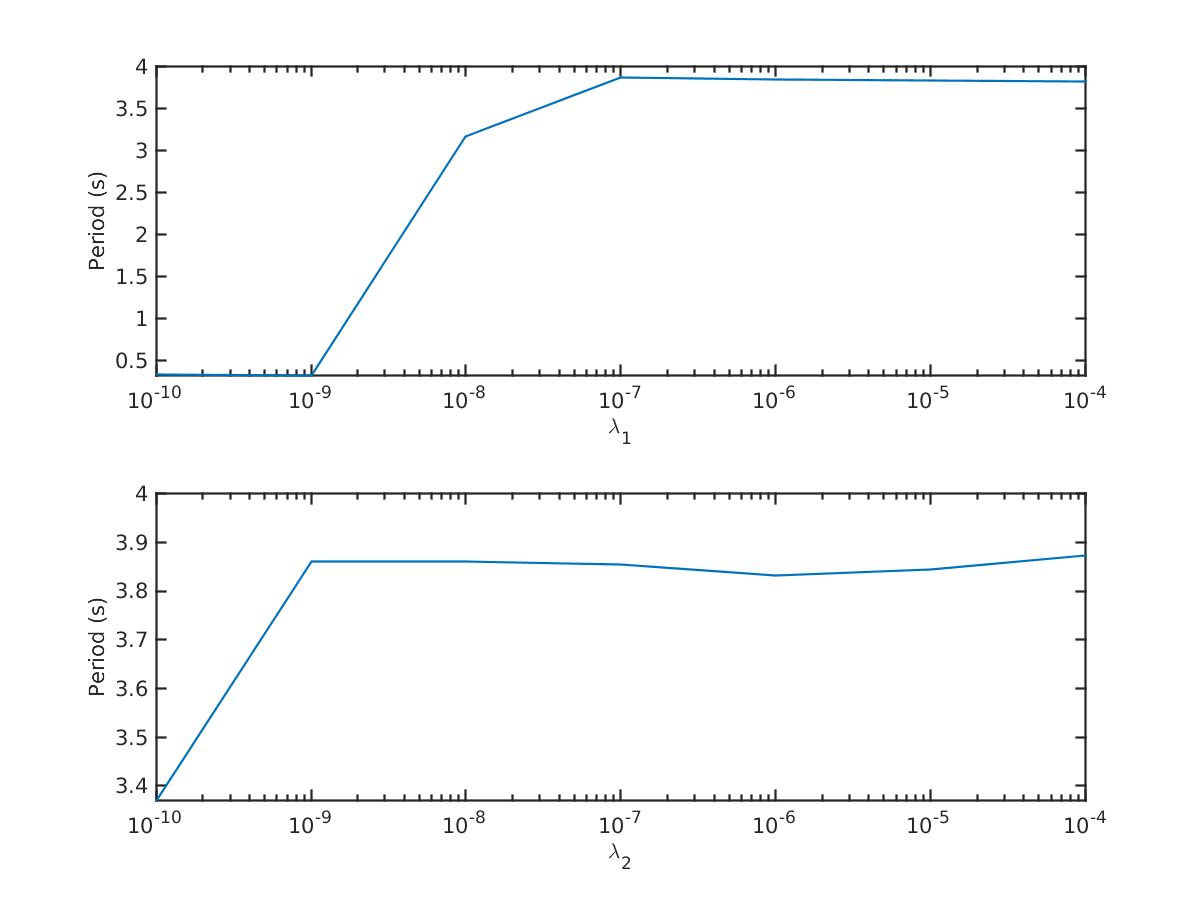}
\caption{Average periods of velocity with different $\lambda_1$ and $\lambda_2$ values, where the wave input has 6-meter significant wave height and 4-second period. Each $\lambda$ is investigated separately, i.e., when we investigate $\lambda_1$, we set $\lambda_2 = 0$.}
\label{average_period_velocity_4sec}
\end{figure}

For these experiments, the optimization time horizon was set as 7 times the wave period, $\Delta t_k = 0.01$, $\gamma = 10^6$, $\delta = 3$, and $\eta = 1$ (meaning that the soft bound equals the hard bound).  We then considered two sets of pairs $(\lambda_1,\lambda_2)$ as follows.
\begin{enumerate}
    \item[(a)] $\lambda_2 = 0$ with $\lambda_1 \in \{10^{-10},10^{-9},\dots, 10^{-4}\}$.
    \item[(b)] $\lambda_1 = 0$ with $\lambda_2 \in \{10^{-10},10^{-9},\dots, 10^{-4}\}$.
\end{enumerate}
Using these pairs of values, we solved the optimization problem and calculated the average realized period of the velocity trajectory.  Figure \ref{average_period_velocity_4sec} shows the calculated average period of the velocity trajectory for a wave input with a period of 4 seconds and a significant wave height of 6 meters.

The figure bears out the suggestion above that the calculated velocity trajectory has a shorter period than the waves if the coefficients are small, and approaches the period of the waves as the coefficients increase. We choose the smallest $\lambda_1$ and $\lambda_2$ values such that the average period of the velocity is within $5\%$ of the wave period.  For this instance, we chose $\lambda_1 = 10^{-7}$ and $\lambda_2 = 10^{-9}$.

\begin{figure}
\centering
\includegraphics[width=10cm]{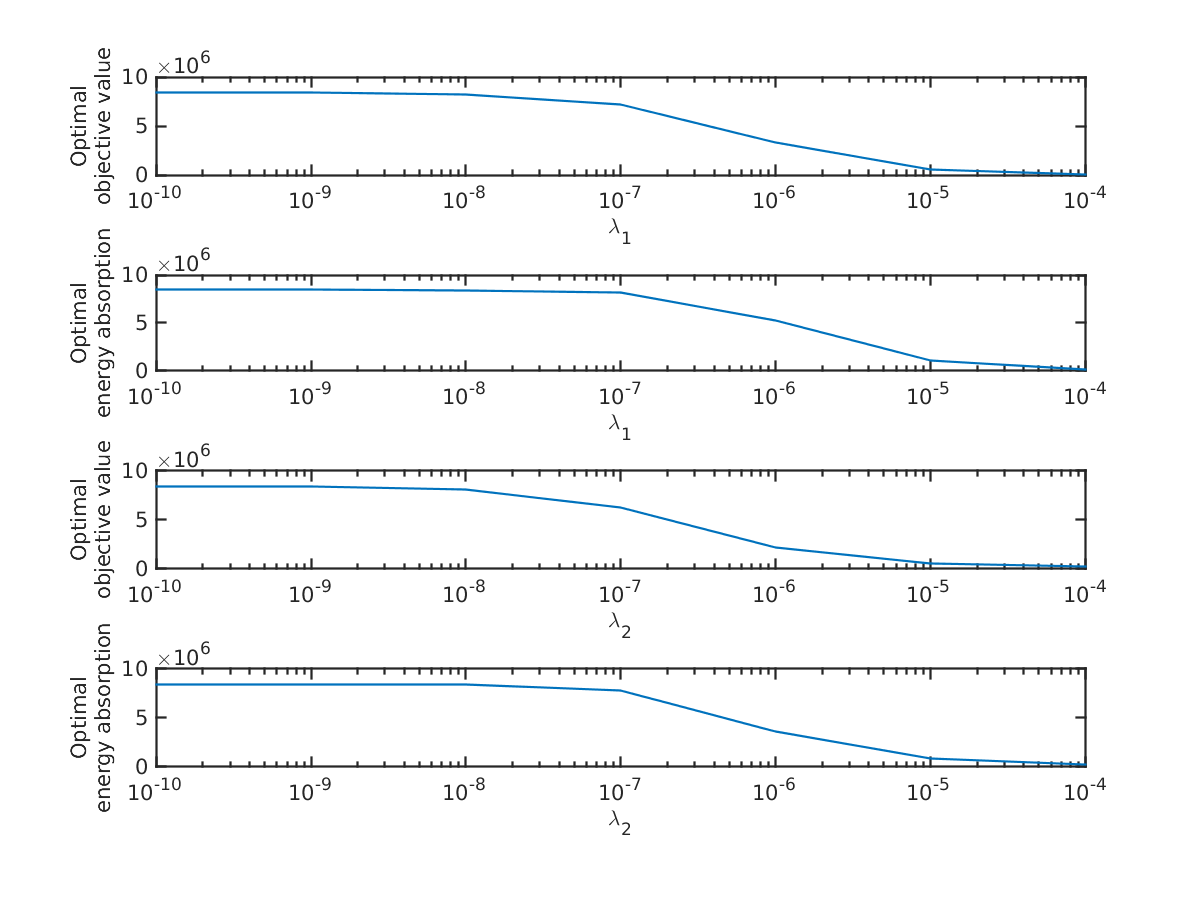}
\caption{Optimal objective values and energy absorption with different $\lambda_1$ and $\lambda_2$ values, where the wave input has 6-meter significant wave height and 4-second period. Each $\lambda$ is investigated separately, i.e., when we investigate $\lambda_1$, we set $\lambda_2 = 0$.}
\label{obj_lambda_4sec}
\end{figure}

As a sanity check, Figure~\ref{obj_lambda_4sec} provides optimal objective values and energy absorption over different $(\lambda_1,\lambda_2)$ values.  The figure shows that with larger cost coefficients, the energy absorption and the objective function go to 0, suggesting that the optimal control trajectory goes to 0, which makes sense for large $(\lambda_1,\lambda_2)$ values.  Overall, these results help us to verify that $\lambda_1 = 10^{-7}$ and $\lambda_2 = 10^{-9}$ are large enough to keep the period of the device in sync with the wave, yet small enough that the device is capable of absorbing energy in the model.

When we conducted similar experiments considering the other wave inputs, we found that slightly different $(\lambda_1,\lambda_2)$ values satisfied our criterion.  This is due to the fact that by changing the period of the wave input, the optimal control sequence has a different period, thus the cost for applying control might require a different coefficient to ensure that the period of the device and the wave are in sync.  Since, in a realistic setting, the parameters $\lambda_1$ and $\lambda_2$ are likely to be fixed (not tunable for different wave conditions), as a final choice of $\lambda_1$ and $\lambda_2$ we choose the largest values among all trials with different wave inputs. This led us to choose $\lambda_1^* = 10^{-6}$ and $\lambda_2^* = 10^{-6}$ for our remaining experiments.

As a further sanity check, as well as to illustrate the sensitivity of the optimal solution with respect to the values of $\lambda_1$ and $\lambda_2$, we evaluate the error that arises from using incorrect cost coefficients in the model. 
In particular, suppose that the true cost coefficients are $\lambda^* = (\lambda_1^*,\lambda_2^*) = (10^{-6},10^{-6})$. 
We then define the following quantities:
\begin{align*}
  f_{\lambda^*}(x(\lambda)) &: \text{objective function evaluated with $\lambda^*$ using the solution obtained} \\ & \text{by solving the problem with cost coefficients $\lambda = (\lambda_1,\lambda_2)$} \\
  f_{\lambda^*}(x(\lambda^*)) &: \text{objective function evaluated with $\lambda^*$ using the solution obtained} \\ & \text{by solving the problem with cost coefficients $\lambda^* =  (\lambda_1^*,\lambda_2^*)$} 
\end{align*}

\begin{figure}
\centering
\includegraphics[width=10cm]{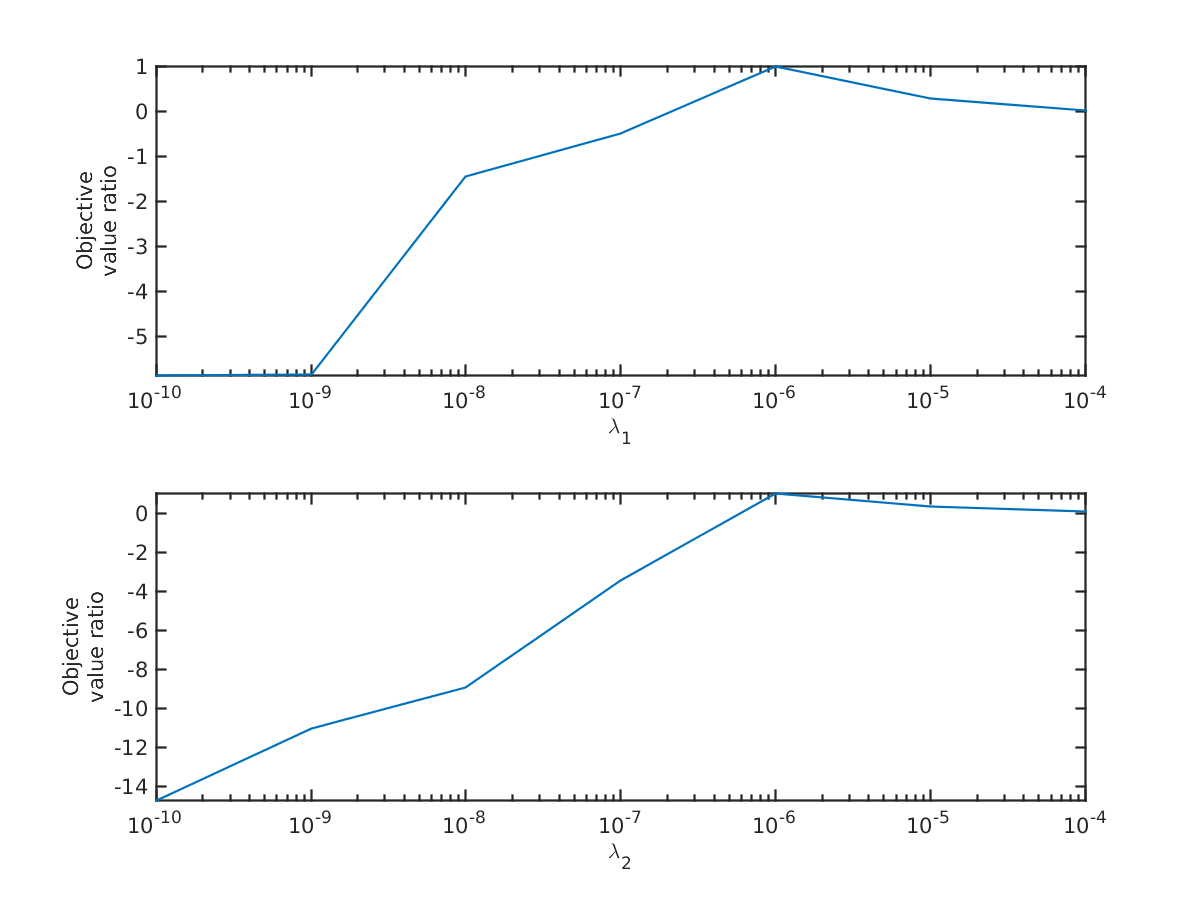}
\caption{Assuming that the true cost coefficients are $\lambda_1^*=10^{-6}$ and $\lambda_2^*=10^{-6}$, objective value ratios with different $\lambda_1$ and $\lambda_2$ values, where the wave input has a significant wave height of 6 meters and a wave period of 4 seconds.}
\label{sensivity_lambda}
\end{figure}

In Figure~\ref{sensivity_lambda}, we plot the values
\begin{equation*}
    \text{ratio} = \frac{f_{\lambda^*}(x(\lambda))}{f_{\lambda^*}(x(\lambda^*))}
\end{equation*}
as a function of $\lambda_1$ with $\lambda_2 = \lambda_2^*$ (top plot) and also as a function of $\lambda_2$ with $\lambda_1 = \lambda_1^*$ (bottom plot).  When $\lambda_1 = \lambda_1^*$ (resp.~$\lambda_2=\lambda_2^*$), the ratio is equal to~1 in the top plot (resp.~bottom plot).  Otherwise, the ratio shows the objective value lost by solving the problem with the ``wrong'' cost coefficients.  We see from this experiment that it is arguably much worse to underestimate the cost coefficients than to overestimate them.  In particular, if the problem is solved with a cost coefficient that is too large, then the plotted ratio is less than 1, but at least still positive.  On the other hand, if the problem is solved with a cost coefficient that is too small, then the plotted ratio can be negative.  These experiments support our inclusion of these costs in the objective function, and support our idea of choosing the largest coefficients among all trials.

\subsection{Choosing proper $\eta$ and $\rho$ values promoting reduced damage to the device}
\label{choosing-safety-parameters}
In order to decide on appropriate values for $\eta$ (the soft bound ratio) and $\rho$ (the penalty term for exceeding the bound), we model and solve the optimization problem \eqref{complete_obj}--\eqref{nonneg} 
under multiple values of $\eta$ and $\rho$. In particular, we set 
 the problem horizon as one full wave period, $\Delta t_k = 0.01$, $\gamma = 10^6$, and $\delta = 3$, and $\lambda$ values set as in the previous section. We tested all combinations of $\eta$ values in $\{\frac{1}{50},\frac{2}{50}, \dots,\frac{49}{50}, \frac{50}{50}\}$ and $\rho$ values in $\{\frac{1}{50},\frac{2}{50}, \dots,\frac{49}{50}, \frac{50}{50}\}$, resulting in  a total of $50 \times 50 = 2500$ instances. The whole procedure was relatively fast since we limit the horizon to one full wave period, and since IPOPT can solve the model \eqref{complete_obj}--\eqref{nonneg} quickly.
It is worth mentioning that all the instances are successfully solved to local optimality. Note that each of the instances is solved multiple times with different initial points, and in all of these cases we get the same optimal solution, suggesting that the local optimum may be globally optimal.

\begin{figure}
\centering
\includegraphics[width=10cm]{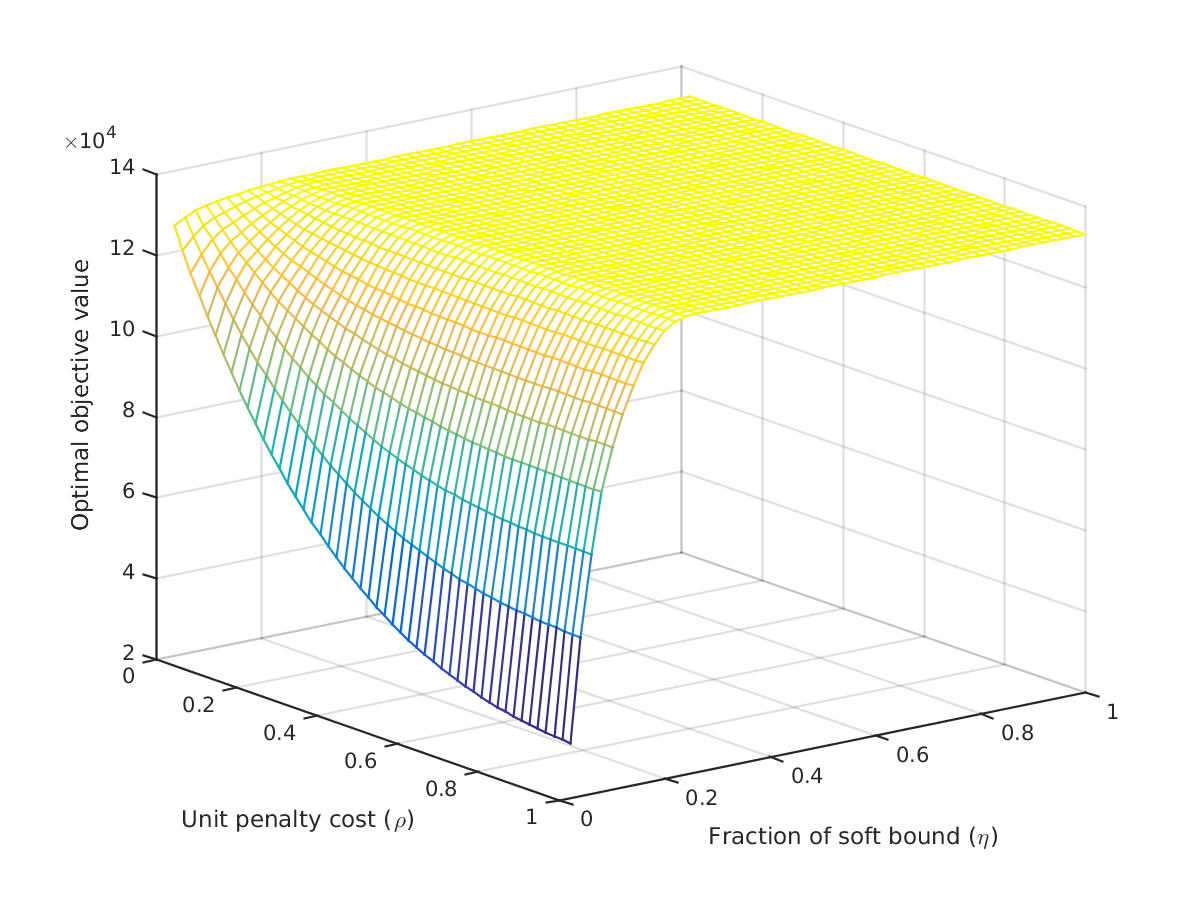}
\caption{Optimal objective values with different $\eta$ and $\rho$ values, where the wave input has 6-meter significant wave height and 4-second period.}
\label{obj_safety_4sec}
\end{figure}

Figure \ref{obj_safety_4sec} plots the optimal objective function value under different safety parameters. As can be seen from the figure, the optimal objective value decreases as $\rho$ increases and as $\eta$ decreases. Moreover, the optimal objective function value is flat for $\eta$ greater than $\approx 0.2$, since below that range, the optimal control trajectory is already bounded away from $\gamma$. This effect will be demonstrated in more detail in the next section. Also, it can be said that the effect of $\eta$ is more direct, and that it inherently affects the impact of $\rho$ on the objective function.

Of course, the choice of the safety parameters primarily depends on the application and is subject to the modeler's preferences. We will further explore the impact of $\eta$ and $\rho$ on the optimal solution in the next section.

\begin{figure}
\centering
\includegraphics[width=10cm]{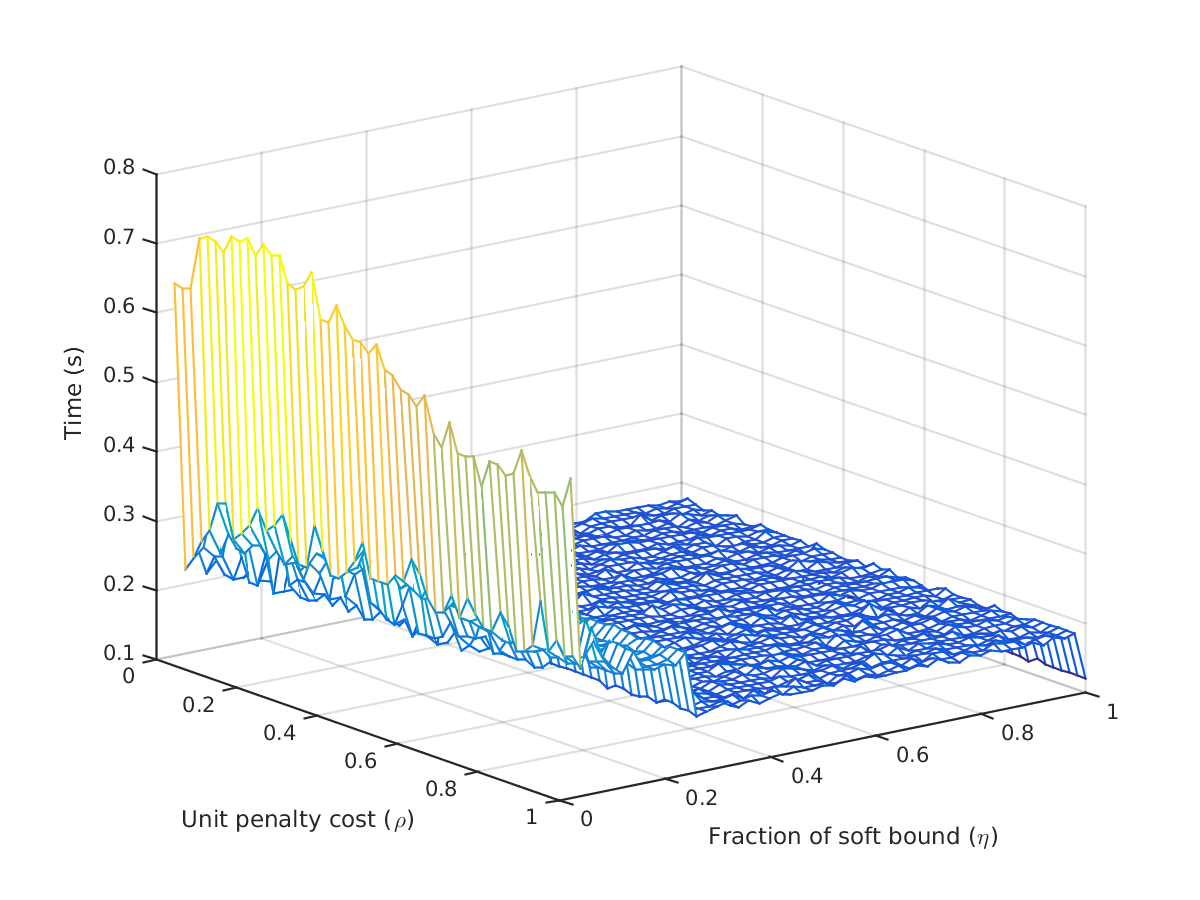}
\caption{CPU times with different $\eta$ and $\rho$ values, where the wave input has 6-meter significant wave height and 4-second period.}
\label{cpu_safety_4sec}
\end{figure}

Figure \ref{cpu_safety_4sec} plots the total CPU time (user CPU + system CPU) returned by the solver IPOPT integrated in the modelling software AMPL. The total CPU time required to solve the problem is relatively small throughout the parameter grid, except when the fraction of soft bound ($\eta$) gets small and becomes an active constraint. In this region, the problem becomes dramatically harder as $\eta \rightarrow 0$, with CPU times roughly quadrupled compared to the other instances.

\subsection{Receding horizon experiments}
\label{receding-horizon-experiments}
In this section, we use the receding horizon approach described in Section \ref{real-time}, using the previously determined $\lambda$ values, initial time $T_0 = 200$, horizon $T = \text{wave period}$, time step $\Delta t = 0.01$, controller update horizon $T_c = \frac{1}{10}\times(\text{wave period})$, and number of receding periods $K = 100$. We also test each of the following $\eta$ and $\rho$ values:
\begin{enumerate}
    \item $\eta =1$ (soft bound equals hard bound),
    \item $\eta =0.1$ and $\rho = 1$,
    \item $\eta =0.1$ and $\rho = 0.1$.
\end{enumerate}

In order to ensure that the control sub-sequence resulting from each problem at the end of a receding period is connected to the next one, we set the initial point of the control trajectory in each receding period as the last entry of the implemented control trajectory obtained from the solution of the previous problem.

\subsubsection{Receding horizon results when $\eta=1$}

\begin{figure}
\centering
\includegraphics[width=10cm]{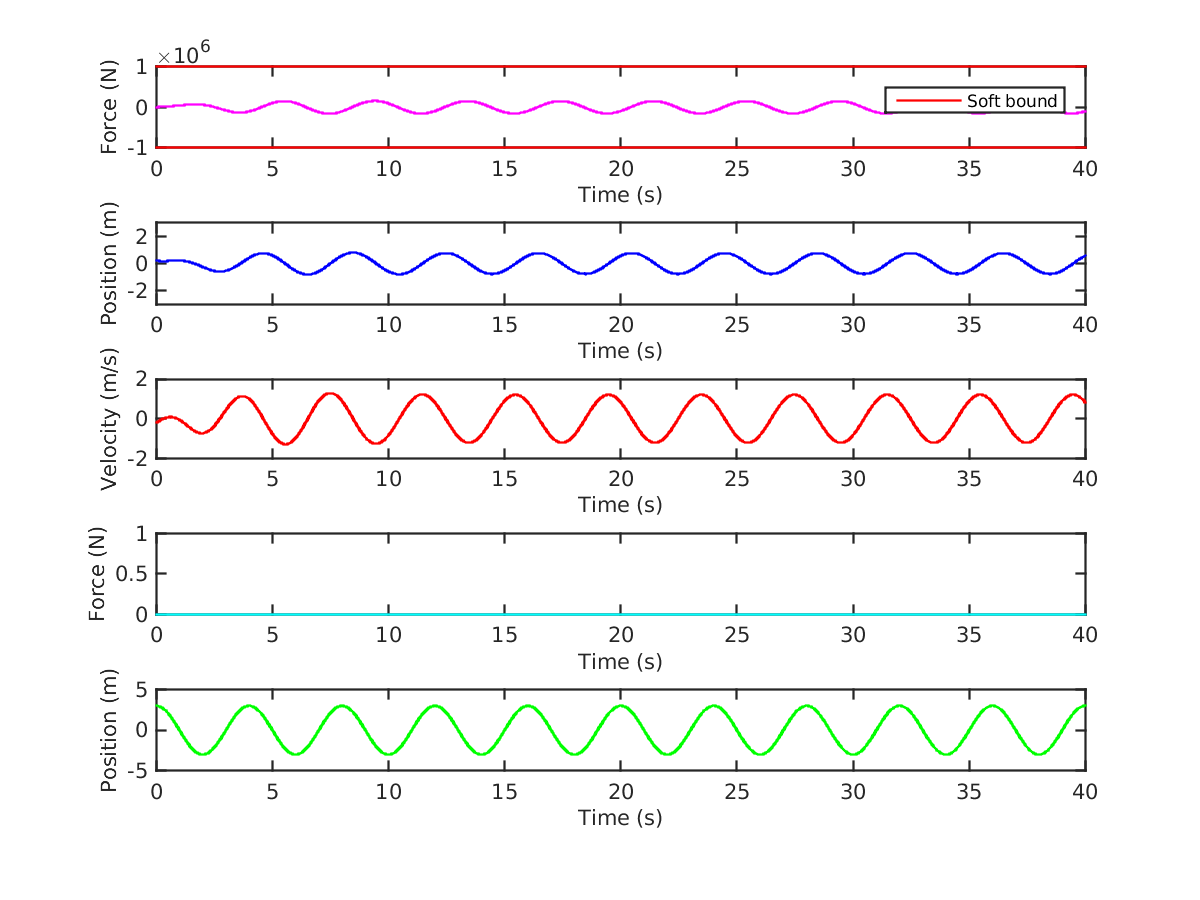}
\caption{Solution of the receding horizon approach with $\eta =1$, where the wave input has 6-meter significant wave height and 4-second period. The subplots show control ($u$), position ($z$), velocity ($\dot{z}$), excess force ($\alpha$), and wave elevation ($w$), respectively.}
\label{soln_real_time_no_safety_4sec_eta_1_rho_1}
\end{figure}

Figure \ref{soln_real_time_no_safety_4sec_eta_1_rho_1} plots the solution of the real time procedure (first plot), and the resulting trajectory (second and third plots), the excess force (fourth plot), and, for reference, the incident wave elevation (fifth plot). As can be seen, the control is smooth and has the same period as the wave input. As discussed in the previous section, the control is bounded away from the hard bound $\gamma$, due to the choice of $\lambda$ parameters. Since the soft bound is chosen as the hard bound, the excess force (i.e., $\alpha_k$, the amount by which the force exceeds the soft bound) is zero. The velocity also has the same period as the wave input, which was guaranteed by our proper choice of $\lambda$ values.

\begin{figure}
\centering
\includegraphics[width=10cm]{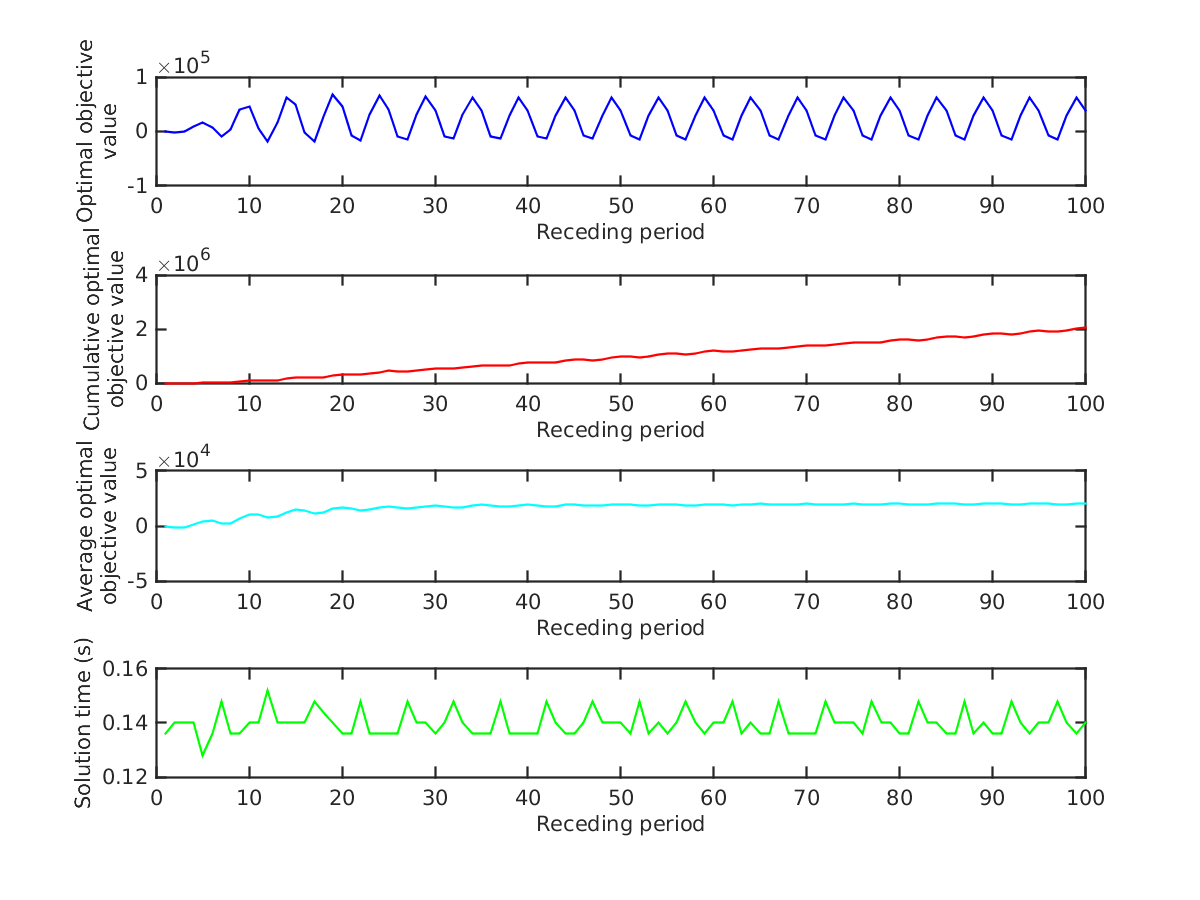}
\caption{Summary of the receding horizon approach with $\eta =1$, where the wave input has 6-meter significant wave height and 4-second period.}
\label{summary_real_time_no_safety_4sec_eta_1_rho_1}
\end{figure}

Figure \ref{summary_real_time_no_safety_4sec_eta_1_rho_1} plots a summary of each problem solved at each receding period. The $x$-axis represents the problem at each receding period, where the total number of receding periods is $K = 100$. Note that at each receding period, only the objective value corresponding to the implemented control policy is calculated. We can see the periodicity in the optimal objective value every 5 receding periods since the controller update horizon is chosen to be $\frac{1}{5}\times\text{(half wave-period)}$ and the generation repeats itself every half wave-period. More importantly, in order to maximize the total energy extracted, the optimal objective value becomes negative for short periods of time. This result supports the remark made by Falnes~\cite{falnes2002ocean}, which states that in optimal control, the control system must be efficient enough to feed back part of the extracted energy. Moreover, we can see that the average optimal objective value becomes flat after a few periods since the effect of the initial point decreases over time. It is also clear that the solution time fluctuates around 0.14 seconds, never exceeding 0.16 or falling below 0.12. Therefore, this specific experiment can be done in real time since the solution time corresponding to each receding horizon is around $0.14s$, which is strictly less than the controller update horizon $\frac{1}{10} 4 = 0.4s$.

\subsubsection{Receding horizon results when $\eta =0.1$ and $\rho = 1$}

\begin{figure}
\centering
\includegraphics[width=10cm]{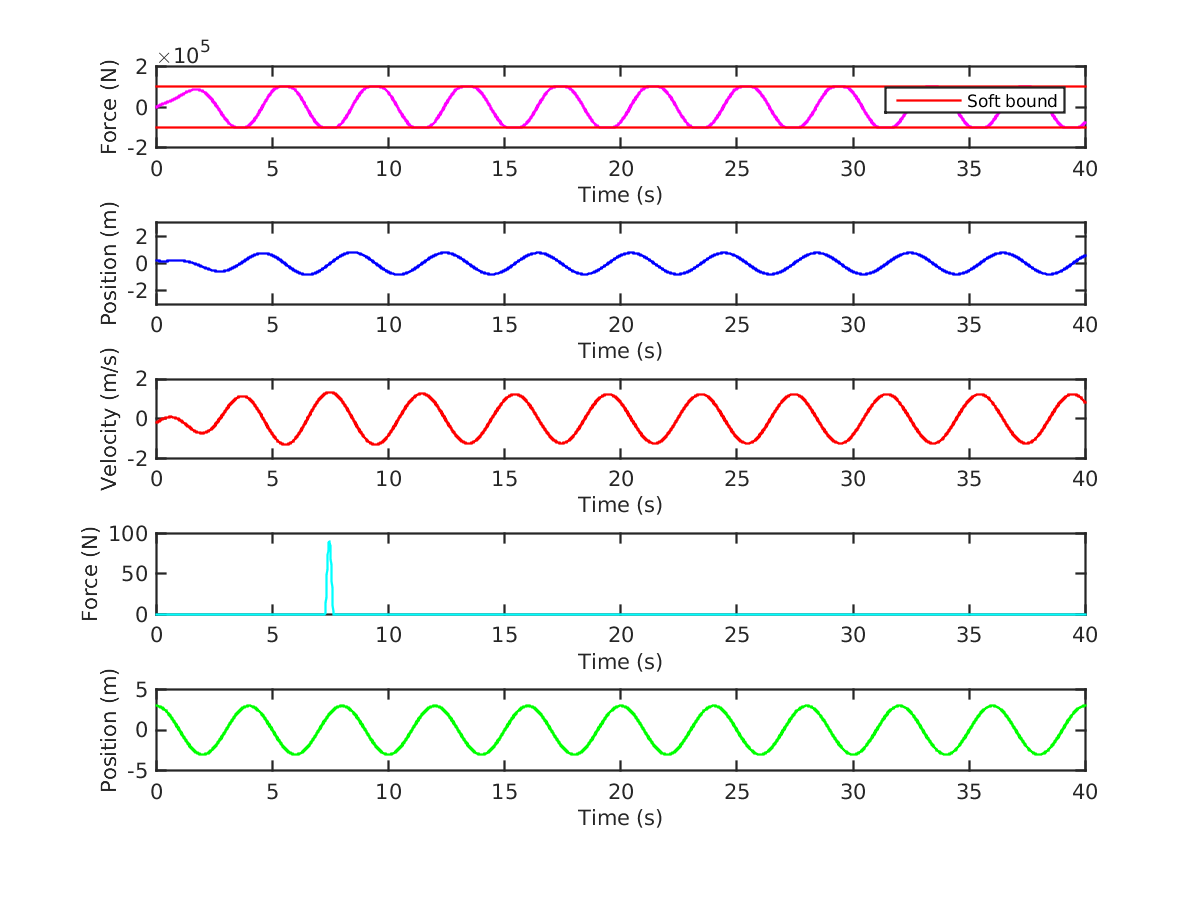}
\caption{Solution of the receding horizon approach with $\eta =0.1$ and $\rho=1$, where the wave input has 6-meter significant wave height and 4-second period. The subplots show control ($u$), position ($z$), velocity ($\dot{z}$), excess force ($\alpha$), and wave elevation ($w$), respectively.}
\label{soln_real_time_no_safety_4sec_eta_01_rho_1}
\end{figure}

Whereas Figure \ref{soln_real_time_no_safety_4sec_eta_1_rho_1} considered the extreme case in which the soft bound equals the hard bound, 
Figure \ref{soln_real_time_no_safety_4sec_eta_01_rho_1} plots another extreme case, in which the soft bound imposes a very tight constraint on the control trajectory ($\eta=0.1$). The choice of $\rho =1$ makes the cost of exceeding the soft bound very large, hence, the control is essentially restricted to the soft bound, with one exception around 7 seconds into the horizon. This violation is very small: the bound of $2\times 10^5$ N is violated by approximately $100$ N. Other than that, the position and velocity are similar to those in Figure \ref{soln_real_time_no_safety_4sec_eta_1_rho_1} except that now the motion is slightly less dampened.

\begin{figure}
\centering
\includegraphics[width=10cm]{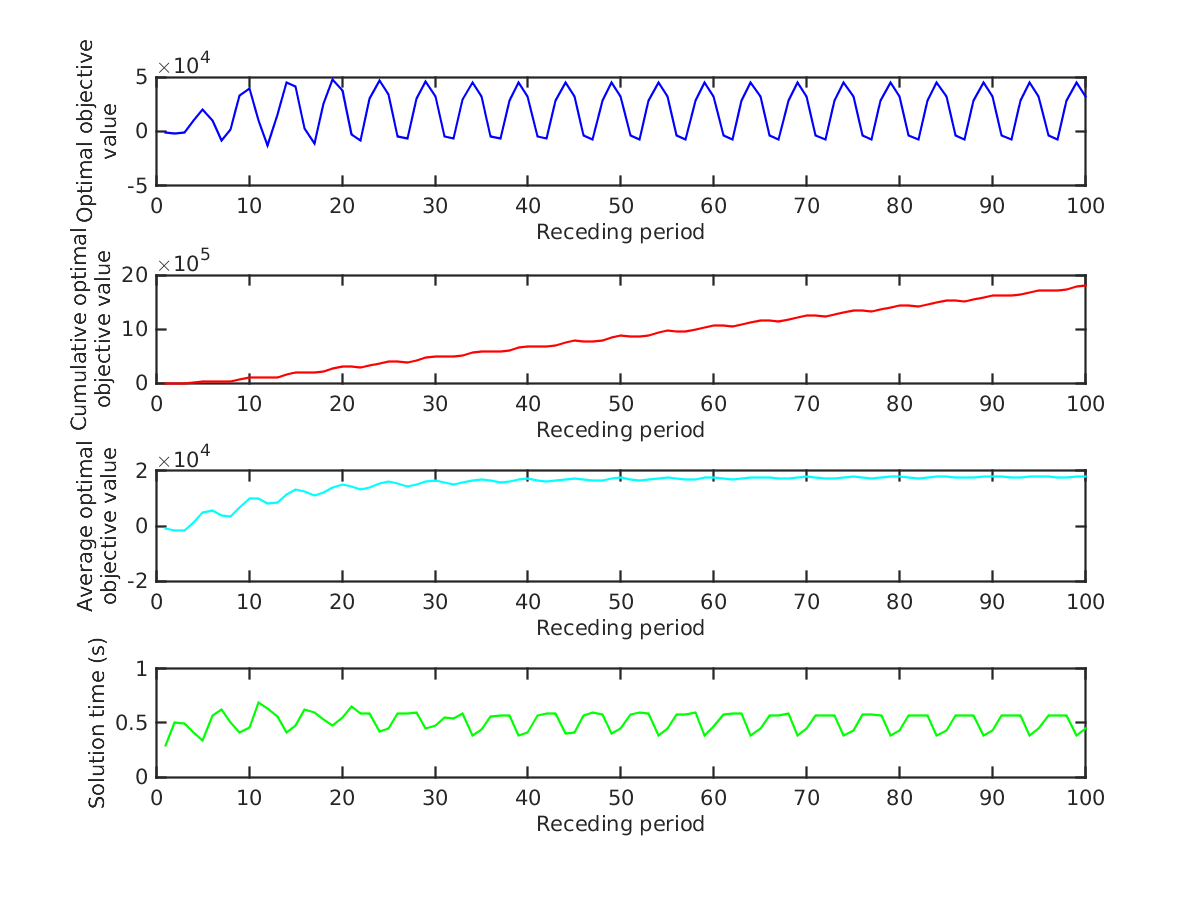}
\caption{Summary of the receding horizon approach with $\eta =0.1$ and $\rho =1$, where the wave input has 6-meter significant wave height and 4-second period.}
\label{summary_real_time_no_safety_4sec_eta_01_rho_1}
\end{figure}

Figure \ref{summary_real_time_no_safety_4sec_eta_01_rho_1} gives us a similar summary, in which now the problem at each receding period is much harder to solve, with an average solution time of $0.5s$ (not a surprise, in light of Figure \ref{cpu_safety_4sec}). Moreover, the objective value has decreased compared to the case in which $\eta=1$, since now we are limiting the amplitude of the control by having a large cost associated with exceeding the soft bound. Unfortunately, this instance cannot be implemented in real time, since solving each problem requires more time than the the controller update horizon. However, one can adjust the controller update horizon appropriately if real-time control is essential. 

\subsubsection{Receding horizon results when $\eta =0.1$ and $\rho = 0.1$}
Now, we consider the same soft bound used in the previous section with a much smaller cost penalty for violating it:  $\rho=0.1$.

\begin{figure}
\centering
\includegraphics[width=10cm]{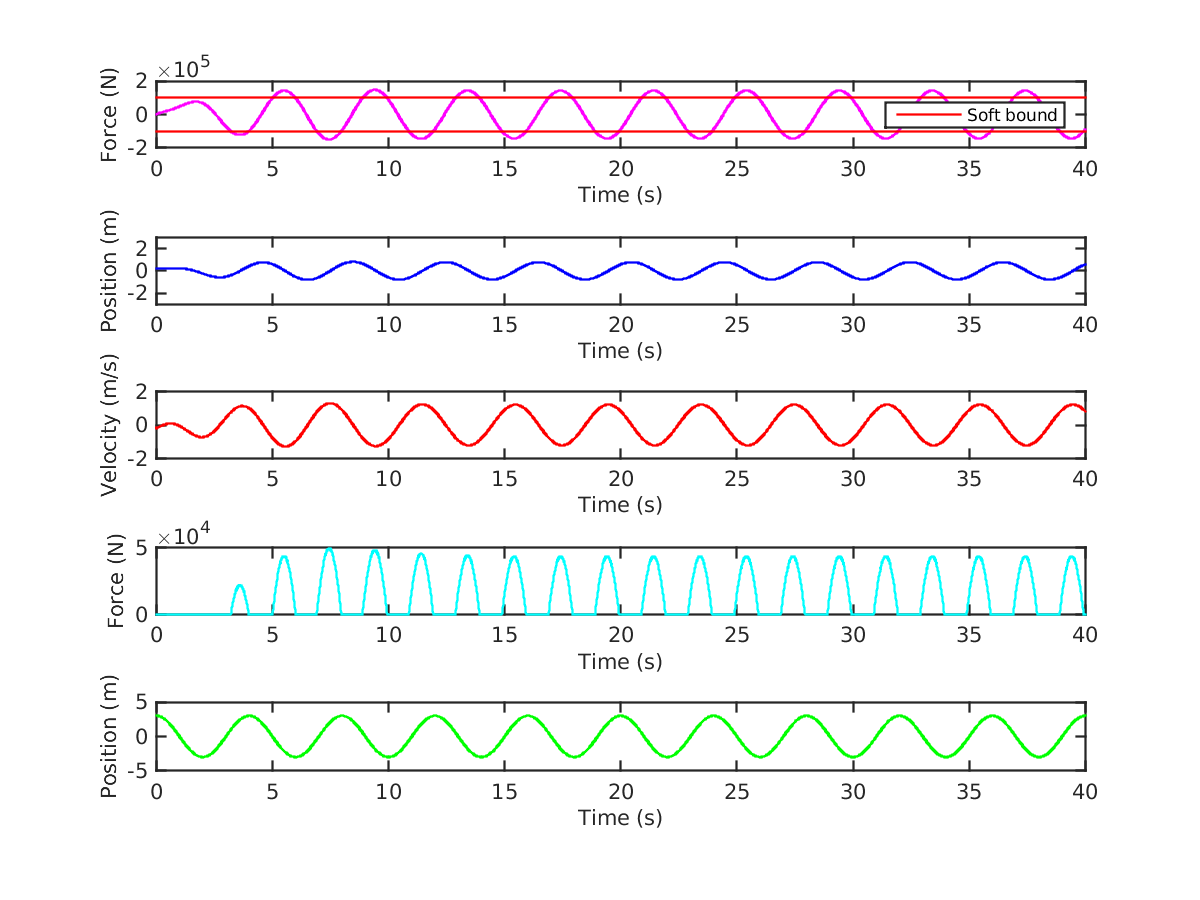}
\caption{Solution of the receding horizon approach with $\eta =0.1$ and $\rho=0.1$, where the wave input has 6-meter significant wave height and 4-second period. The subplots show control ($u$), position ($z$), velocity ($\dot{z}$), excess force ($\alpha$), and wave elevation ($w$), respectively.}
\label{soln_real_time_no_safety_4sec_eta_01_rho_01}
\end{figure}

As  can be observed in Figure \ref{soln_real_time_no_safety_4sec_eta_01_rho_01}, there aren't large violations of the soft bound, which can be seen more clearly in the excess force plot, where it is nonzero around the peak points of the control trajectory. Since we used a non-zero penalty for the violation $(\rho = 0.1)$, this result demonstrates a trade-off between exceeding the soft bound and absorption of energy. In other words, we can observe that the trajectory amplitude here is strictly between those in Figures \ref{soln_real_time_no_safety_4sec_eta_1_rho_1} and \ref{soln_real_time_no_safety_4sec_eta_01_rho_1}.

\begin{figure}
\centering
\includegraphics[width=10cm]{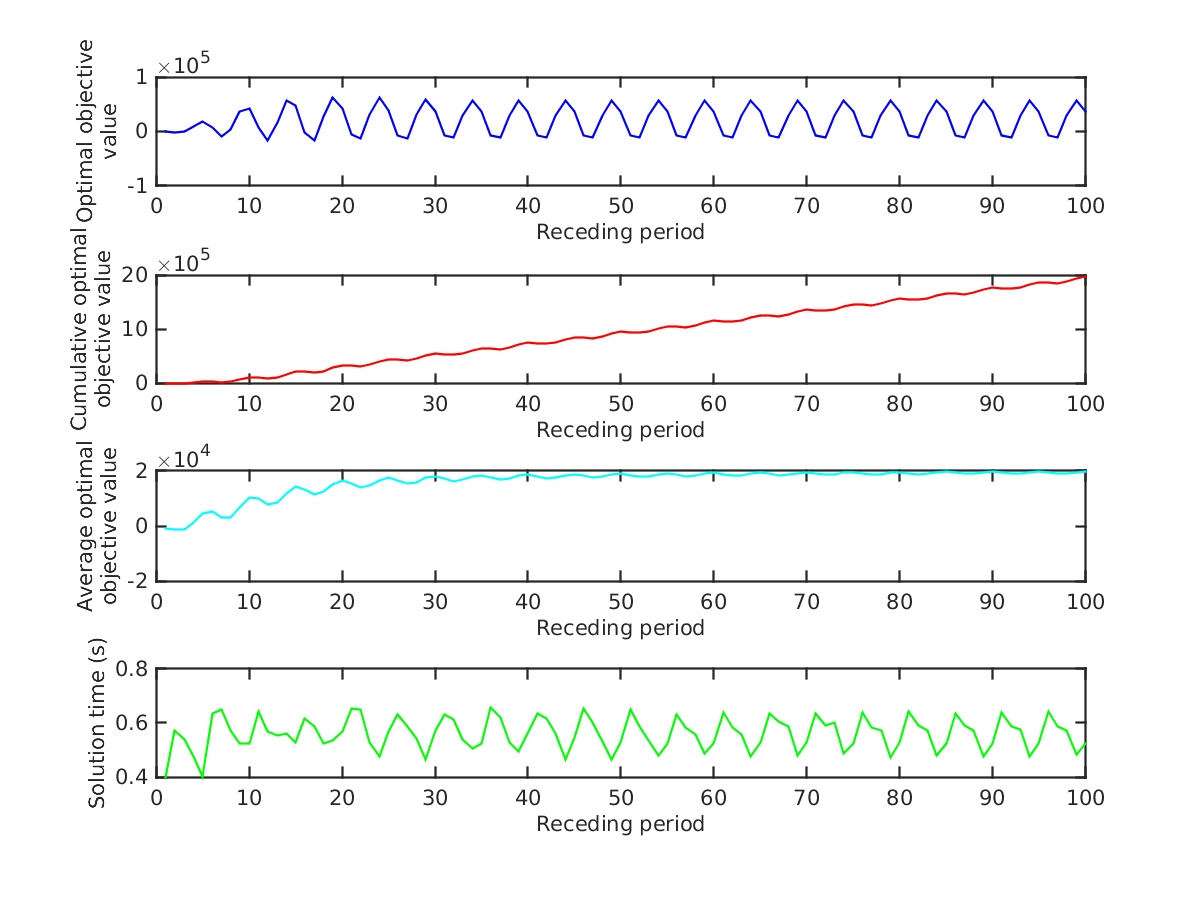}
\caption{Summary of the receding horizon approach with $\eta =0.1$ and $\rho = 0.1$, where the wave input has 6-meter significant wave height and 4-second period.}
\label{summary_real_time_no_safety_4sec_eta_01_rho_01}
\end{figure}

Based on Figure \ref{summary_real_time_no_safety_4sec_eta_01_rho_01}, the computational difficulty of the problem remains similar to that in Figure \ref{summary_real_time_no_safety_4sec_eta_01_rho_1}; however, in this case, the objective is slightly better since we impose less of a penalty for exceeding the soft bound. Additionally, as was expected, the objective value is worse than the case given in Figure \ref{summary_real_time_no_safety_4sec_eta_1_rho_1} since the soft bound was not active in that case. Finally, due to the relatively large solution times, this instance is not real-time solvable.

\section{Conclusion}
\label{conclusion}

In this paper, we propose a generic optimal control model for a single wave-energy converter, considering the device to be a two-body point absorber. Our model makes two main contributions to the literature: (1) It incorporates a realistic term to represent the cost or energy expended in the control force, thus addressing a deficiency of earlier models that led to inaccurate estimates of the actual energy extracted from the device. (2) It imposes both hard and soft bounds on the control force and relative displacement, thus promoting safe operation of the device.

The hydrodynamic coefficients of the device are estimated using simulation via WEC-Sim. We derived an approximate form for the cost of the control force and integrated it into the objective function. The formulation of the optimal control problem resulted in a quadratic program (QP) with a non-convex objective function. 

In our numerical results, we  demonstrate the effect of cost and safety parameters on the optimal energy absorption, underlining the importance of correctly estimating the cost parameters. Furthermore, we illustrate the trade-off between the safe operation of the device and the energy absorbed. Finally, we display the capabilities of the model by utilizing it in a real-time optimization framework.

Our model and analysis assume regular waves. An important extension for future work is to consider irregular waves as the wave input. This extension requires incorporation of a prediction scheme and the update of system coefficients into the real-time optimization framework. Hopefully, the update of the system coefficients can be done in a systematic manner via appropriate transformation, rather than using an MPC-like approach, to further decrease the complexity of the real-time implementation.

Another possible extension is the control of multiple WECs simultaneously, while maintaining the optimality and the quality of the estimation. This problem is rather challenging since it requires a dynamic estimation of the coefficients incorporating the interaction between bodies, and implementing a control policy that can optimize the system as a whole rather than individually.

We also believe that a better form for the cost of the control force may be derived, by utilizing more detailed information of a specific WEC. The incorporation of a more accurate form for the cost of the control force would further reduce the gap between the theory and the practice of WEC control.

\section*{Acknowledgement}

This research was supported in part by NSF grants \#CMMI-1400164 and CCF-1442858. This support is gratefully acknowledged.

\section*{References}

\bibliography{mybibfile}

\begin{appendices}

\section{Derivation of functional form for the cost of the control force}
\label{derivation}
\begin{figure}
\centering
\includegraphics[width=10cm]{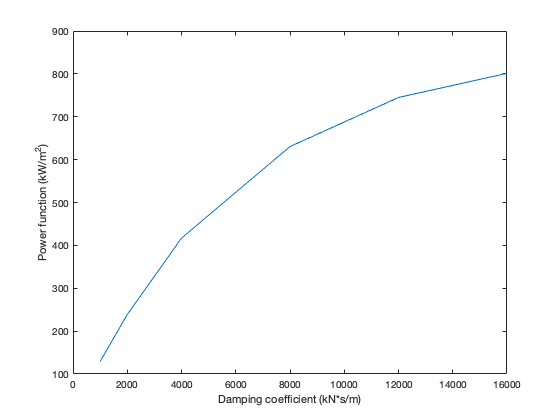}
\caption{Power function under different damping coefficients}
\label{powerfunc}
\end{figure}

The functional form of the cost of the control force is derived by considering a damper as a controller. For simplicity, in the derivation, we will use the notation that is introduced in \cite{klamo2007effects} (rather than attempting to bring their notation in line with ours), and explain the notation as it is introduced. In \cite{klamo2007effects}, the author derives the magnetic damping term as
\begin{equation}
    \label{magnetic-damping}
    b_{mag} = f(w,R,\dots) i_{sup}^2,
\end{equation}
where $i_{sup}$ is the current supplied to the system and $f(w,R,\dots)$ is a constant depending on the magnetic system parameters. The power required to generate current $i_{sup}$ can be found (under an assumption of DC power) through
\begin{align*}
    P &= VI, 
    \intertext{which can be rewritten using Ohm's law, $V = IR$, to obtain}
    P &= I^2 R \\
    I^2 &= \frac{P}{R},
\end{align*}
where $P$ is the power, $V$ is the voltage, $I$ is the current and $R$ is the resistance. Plugging the last equation into \eqref{magnetic-damping} gives
\begin{align}
    b_{mag} &= \Tilde{f}(w,R,\dots) P, \nonumber
    \intertext{where $\Tilde{f} = \frac{f}{R}$. This gives us the relation}
    \label{power-damping}
    P &= \frac{b_{mag}}{\Tilde{f}}.
\end{align}
Ultimately, we want to obtain a functional form such that for some function $g$, we have $P \approx g(F)$, where $F$ represents the applied control force and $P$ represents the power required to generate that force. From now on, our focus will only be on the functional form of $g$, without considering the constant terms. To this end, we will use the form of a damping controller,
\begin{equation}
\label{control-form}
    F = -b \dot{z},
\end{equation}
where $\dot{z}$ represents the velocity and $b$ the damping coefficient, then try to capture the relation between $b$ and $\dot{z}$ under the best possible sea conditions for the two-body device we are considering. By using the force equation, we have the generated power expressed as
\begin{equation*}
    P = -F \dot{z} = b \dot{z}^2.
\end{equation*}
Now, based on the experiments in \cite{muliawan2013analysis}, we have some experimental data for the power capture under different damping coefficients.  Specifically, Figure~\ref{powerfunc} shows several power function values when different damping coefficients are implemented in the control strategy, where the control strategy has the form given in \eqref{control-form}. It is important to note that the power function is wave-normalized for regular waves.

\begin{figure}
\centering
\includegraphics[width=10cm]{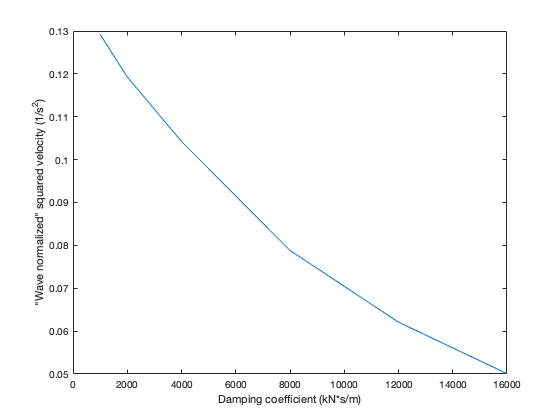}
\caption{Average squared velocity under different damping coefficients}
\label{v2}
\end{figure}

By dividing the normalized power by the damping coefficient, we obtain the values shown in Figure~\ref{v2}.

\begin{figure}
\centering
\includegraphics[width=10cm]{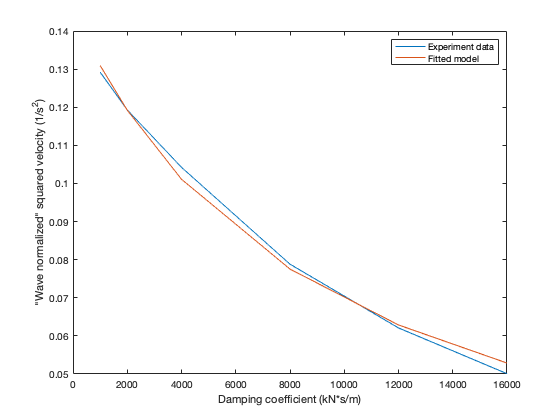}
\caption{The fitted model for the relation between average squared velocity and the damping coefficient}
\label{model-fit}
\end{figure}

The next step is to fit a functional relationship between $b$ and the average $\dot{z}^2$, but before doing that, let us observe some properties of the relation between them:
\begin{enumerate}
    \item When $b \shortuparrow$, $\dot{z}^2 \shortdownarrow$.
    \item When $b \rightarrow \infty$, $\dot{z}^2 \rightarrow 0$.
    \item When $b \rightarrow 0$, $\dot{z}^2 \rightarrow$ some upper limit.
    \item The damping coefficient $b$ should increase faster than the decrease in $ \dot{z}^2$, so that the multiplication $b\dot{z}^2$ increases as shown in Figure \ref{powerfunc}.
\end{enumerate}
These properties hold for certain exponential, logarithmic, and polynomial functions. Therefore, we used each of these three families of nonlinear regression models, fitting each to our data. We chose the model with the largest $R^2$ value; that model is:
\begin{equation*}
    \dot{z}^2 = \frac{1330.2}{b + 9158.7},
\end{equation*}
which is graphically shown in Figure \ref{model-fit}.

Now, we can further derive
\begin{align*}
    \dot{z}^2 &\approx \mathcal{O}\left(\frac{1}{b}\right),
    \intertext{and since we know that $b>0$, we have}
    |\dot{z}| &\approx \mathcal{O}\left(\frac{1}{\sqrt{b}}\right).
    \intertext{By combining the relations in \eqref{control-form} and \eqref{power-damping}, we end up with}
    |F| = b |\dot{z}| &\approx \mathcal{O}(\sqrt{b}) \approx \mathcal{O}(\sqrt{P})\\
    P &\approx \mathcal{O}(F^2).
\end{align*}
The last statement,  $P \approx \mathcal{O}(F^2)$, simply describes the resulting functional form for the cost of applying the control force. In our notation, $F$ is represented by $u$, and $P$ represents the power consumed to generate a force $u$.

\end{appendices}

\end{document}